\documentclass[12pt]{amsart}
\usepackage{amsthm,amssymb,amsmath,amscd}

\newtheorem{theorem}{Theorem}[section]

\theoremstyle{definition}

\newtheorem{example}[theorem]{Example}
\newtheorem{remark}[theorem]{Remark}

\title[Extended counterpoint]{Extended counterpoint symmetries and continuous counterpoint}
\date{19/02/2015}

\author{Octavio A. Agust\'in-Aquino}
\address{Universidad de la Ca\~nada\\
San Antonio Nanahuatipan Km 1.7 s/n. Paraje Titlacuatitla\\
Teotitlan de Flores Mag\'on, Oaxaca, M\'exico, C.P. 68540.}
\email{octavioalberto@unca.edu.mx}

\begin{document}
\begin{abstract}
A counterpoint theory for the whole continuum of the octave is obtained from Mazzola's model via extended counterpoint symmetries,
and some of its properties are discussed.
\end{abstract}
\maketitle

\section{Introduction}

Mazzola's model for first species counterpoint is interesting because it
predicts the rules of Fux's theory (in particular, the forbidden parallel fifths)
reasonably well. It is also generalizable to microtonal equally tempered scales of even cardinality, and offers alternative
understandings of consonance and dissonance distinct from the one explored extensively in Europe.
In this paper we take some steps towards an effective \emph{extension} of the whole model from a
microtonal equally tempered scale into another, and not just of the mere consonances and
dissonances, as it was done by the author in his doctoral dissertation \cite{oA11}.

First, we provide a definition
of an \emph{extended} counterpoint symmetry that preserves the characteristics of the counterpoint of one
scale in the refined one. Then, we see that the progressive granulation of a specific
example suggest an infinite counterpoint with a continuous polarity, different from the one that Mazzola himself
proposed; a comparison of both alternatives calls for a deeper examination of the meaning of
counterpoint extended to the full continuum of frequencies within the octave.

We must warn the reader that just a minimum exposition of Mazzola's counterpoint model
is done, and hence we refer to his treatise \emph{The Topos of Music} \cite{gM02} (whose notation we use here) and an upcoming
comprehensive reference \cite{AJM15} for further details.

\section{Some definitions and notations}

Let $R$ be a finite ring of cardinality $2k$. A subset $S$ of $R$ of  such that $|S|=k$ is a \emph{dichotomy}. It is often denoted by $(S/\complement S)$
to make the complement explicit. The group
\[
 \overrightarrow{GL}(R) = R\ltimes R^{\times} = \{e^{u} v:u\in R,v\in R^{\times}\}
\]
is called the \emph{affine group} of $R$, its members are the \emph{affine symmetries}. It acts on $R$ by
\[
e^{u}v(x) = vx+u;
\]
this action is extended to subsets in a pointwise manner. A dichotomy $S$ is called \emph{self-complementary}
if there exists an affine symmetry $p$ (its \emph{quasipolarity}) such that $p(S) = \complement S$. A self-complementary dichotomy is
\emph{strong} if its quasipolarity $p$ is unique, in which case $p$ is called its \emph{polarity}.

Of particular interest are the strong dichotomies of $\mathbb{Z}_{2k}$, since this ring models very well the
equitempered $2k$-tone scales modulo octave and Mazzola discovered that the set of classical consonances is a strong
dichotomy. For counterpoint, the self-complementary dichotomies of the dual numbers
\[
 \mathbb{Z}_{2k}[\epsilon] = \{a+\epsilon.b: a,b\in \mathbb{Z}_{2k},\epsilon^{2}=0\}
\]
are even more interesting, since they are used in Mazzola's counterpoint model as counterpoint intervals. More
specifically, given a counterpoint interval $a+\epsilon. b$,
$a$ represents the cantus firmus, and $b$ the interval between $a$ and the discantus, and from every strong dichotomy $(K/D)$
with polarity $p=e^{u}.v$ in $\mathbb{Z}_{2k}$ we can obtain the \emph{induced interval dichotomy}
\[
(K[\epsilon]/D[\epsilon]) = \{x+\epsilon.k:x\in\mathbb{Z}_{2k},k\in K\}
\]
in $\mathbb{Z}_{2k}[\epsilon]$. It is easily proved that, for every cantus firmus, there exists a quasipolarity $q_{x}[\epsilon]$
that leaves its \emph{tangent space} $x+\epsilon. K$ invariant.

A symmetry $g\in\overrightarrow{GL}(\mathbb{Z}_{2k}[\epsilon])$ is a \emph{counterpoint symmetry} of
the consonant interval $\xi=x+\epsilon.k\in K[\epsilon]$ if
\begin{enumerate}
\item the interval $\xi$ belongs to $g(D[\epsilon])$,
\item it commutes with the quasipolarity $q_{x}[\epsilon]$,
\item the set $g(K[\epsilon])\cap K[\epsilon]$ is of maximal cardinality among those obtained with symmetries
that satisfy the previous two conditions.
\end{enumerate}

Given a counterpoint symmetry $g$ for a consonant interval $\xi$, the members of the set $g(K[\epsilon])\cap K[\epsilon]$ are
its \emph{admitted successors}; they represent the rules of counterpoint in Mazzola's model. It must also be noted
that it can be proved that the admitted successors only need to be calculated for intervals of the form $0+\epsilon.k$, and then suitably transposed
for the remaining intervals.

\section{Extending counterpoint symmetries}

Let $(X_{n}/Y_{n})$ be a strong dichotomy in $\mathbb{Z}_{n}$ where
\[
 g_{1}=e^{\epsilon.t_{1}} (u_{1}+\epsilon.u_{1}v_{1}):\mathbb{Z}_{n}[\epsilon]\to \mathbb{Z}_{n}[\epsilon]
\]
is a contrapuntal symmetry for the consonant interval $\epsilon.y\in X_{n}[\epsilon]$, with
$p_{n}=e^{r_{1}} w_{1}$ the polarity of $(X_{n}/Y_{n})$. This means that if $s\in X_{n}$ and $p_{n}[\epsilon]=e^{\epsilon.r_{1}} w_{1}$ is the induced quasipolarity then
\[
 t_{1} = y-u_{1}p_{n}(s)\quad\text{and}\quad p_{n}[\epsilon](\epsilon.t_{1})=g_{1}(\epsilon.r_{1}),
\]
as it is proved in \cite[p. 652]{gM02}. If $a:X_{n}\hookrightarrow X_{an}:x\mapsto ax$ is an embedding
of dichotomies, then
\[
 p_{an}\circ a = a\circ p_{n}
\]
(where $p_{an}=e^{r_{2}} w_{2}$ is the polarity of $(X_{an}/Y_{an})$) and, evidently,
\[
 p_{an}[\epsilon]\circ a = a \circ p_{n}[\epsilon].
\]
In particular, $ar_{1} = r_{2}$.

Suppose there is a symmetry
\[
g_{2}=e^{\epsilon.t_{2}} (u_{2}+\epsilon.u_{2}v_{2}):\mathbb{Z}_{an}[\epsilon]\to \mathbb{Z}_{an}[\epsilon]
\]
such that $a\circ g_{1} = g_{2}\circ a$, then 
\[
 t_{2}=at_{1}\quad\text{and} \quad au_{2} = au_{1}.
\]

From this we deduce
\begin{align*}
 t_{2} = at_{1} &= ay - au_{1}p_{1}(s)\\
 &= ay - u_{2}ap_{n}(s)\\
 &= ay - u_{2} p_{an}(as)
\end{align*}
where $as\in X_{an}$, and
\begin{multline*}
 p_{an}[\epsilon](\epsilon.t_{2}) =p_{an}[\epsilon](\epsilon.at_{1})\\
 =ap_{n}[\epsilon](\epsilon.t_{1}) = ag_{1}(\epsilon.r_{1})\\
  = g_{2}(\epsilon.ar_{1})= g_{2}(\epsilon.r_{2}).
\end{multline*}

This means that $g_{2}$ is almost a contrapuntal symmetry for $\epsilon.ay$, except for
the maximization of the intersection $g_{2}X_{an}[\epsilon]\cap X_{an}[\epsilon]$.
Now we can define a \emph{extended counterpoint symmetry with respect the embedding $a$}
as a symmetry $g_{2}\in \overrightarrow{GL}(\mathbb{Z}_{an}[\epsilon])$ that satisfy
\begin{enumerate}
\item $a\circ g_{1} = g_{2}\circ a$ with $g_{1}$ a (extended or not) contrapuntal symmetry for $\epsilon.y$, and
\item $g_{2}X_{an}[\epsilon]\cap X_{an}[\epsilon]$ has the maximum cardinality among the symmetries with the above property.
\end{enumerate}

Note that extended counterpoint symmetries preserve the admitted successors
of $\epsilon.y\in \mathbb{Z}_{n}[\epsilon]$,
since otherwise the restriction $g_{2}|_{\mathbb{Z}_{n}[\epsilon]}$ of a extended counterpoint symmetry would be a symmetry
such that the intersection $g_{2}|_{\mathbb{Z}_{n}[\epsilon]}X_{n}[\epsilon]\cap X_{n}[\epsilon]$ is bigger than the
corresponding intersection for any counterpoint symmetry. This is a contradiction.

\begin{remark}
In particular, extended counterpoint symmetries always exist in the case of the embedding
$2:\mathbb{Z}_{n}\to \mathbb{Z}_{2n}$, because all the elements of $GL(\mathbb{Z}_{n})$
are coprime with $2$. Thus, for any $\epsilon.y\in \lim_{k\to \infty} X_{2^{k}\cdot n}[\epsilon]$, there exist
a extended contrapuntal symmetry in the limit $\lim_{k\to\infty} \mathbb{Z}_{2^{k}\cdot n}[\epsilon]$
which is the limit of extended counterpoint symmetries.
\end{remark}

\begin{example}
Let $X_{6}=\{0,2,3\}\subseteq \mathbb{Z}_{6}$. The consonant interval $\epsilon.2\in \mathbb{Z}_{6}[\epsilon]$
has $e^{\epsilon.3} (1+\epsilon.3)$ as its only counterpoint symmetry and $15$ admitted
successors. The extended counterpoint symmetries 
of $\epsilon.4\in X_{12}=\{0,1,4,5,6,9\}\subseteq \mathbb{Z}_{12}$ with respect to the embedding $2$ are
$e^{\epsilon.6}.(1+\epsilon.6)$ and $e^{\epsilon.6}.(7+\epsilon.6)$. The number of extended
admitted successors is $48$.
\end{example}

\section{A more detailed example}

In Example 4.11 of \cite{oA11}, it is shown that there exists a strong dichotomy in $\mathbb{Z}_{24}$ that
can be extended progressively (via the embedding Lemma 4.5 of \cite{oA11}) towards a dense
dichotomy in $S^{1}$ with polarity $x\mapsto xe^{i\pi}$, which is the antipodal map. Analogously, the dichotomy
\[
U_{0}=\{0,1,3,\ldots,7,10\}
\]
in $\mathbb{Z}_{16}$ can be completed in each step using the dichotomy
\[
 V_{i}=\{0,\ldots,|U_{i}|-1\},
\]
so we have the inductive definition
\[
U_{i+1}= 2U_{i}\cup (2V_{i}+1), \quad i\geq 1,
\]
which 
is a strong dichotomy of $\mathbb{Z}_{2^{4+i}}$, in each case with polarity $e^{2^{3+i}}$. Note that 
 the injective limit of the $U_{i}$ in $S^{1}$ is dense in one hemisphere.

The standard counterpoint symmetries for 
$U_{0}$ and successively extended counterpoint symmetries for $\mathbb{Z}_{512}$ are listed
in Table \ref{T:Sim}. With ``successively extended'' we mean that they are those who commute with
the extended counterpoint symmetries of $\mathbb{Z}_{256}$, which in turn commute with those of
$\mathbb{Z}_{128}$, and so on down to $\mathbb{Z}_{16}$. In most cases the linear part is
$-1$, and in fact it is remarkable that all of them have no dual component.

\begin{table}
\begin{tabular}{|c|c|c|c|c|}
\hline
Interval & Symmetries & $|gX[\epsilon]\cap X[\epsilon]|$ & Extended & $|gX[\epsilon]\cap X[\epsilon]|$\\
&for $\mathbb{Z}_{16}$ & & symmetries &\\
& & & for $\mathbb{Z}_{512}$ &\\
\hline
$0$ & $e^{\epsilon 5}3$& $96$ &&\\
&$e^{\epsilon 6}13$ &&&\\
&$e^{\epsilon 11}15$ &&$e^{\epsilon 352}511$&$82432$\\
\hline
$1$ & $e^{\epsilon 10}15$ & $112$ &$e^{\epsilon 320}511$&$98816$\\
\hline
$3$ & $e^{\epsilon 2}5$ & $96$&&\\
& $e^{\epsilon 9}11$ & &&\\
& $e^{\epsilon 11}15$&&$e^{\epsilon 352}511$&$82432$\\
\hline
$4$ & $7$ & $112$&$7$&$75264$\\
&&&$439$&\\
\hline
$5$ & $e^{\epsilon 1}3$ & $96$&&\\
&$e^{\epsilon 6}13$&&&\\
&$e^{\epsilon 7}15$&&$e^{\epsilon 244}511$&124416\\
\hline
$6$ & $e^{\epsilon 3}13$ & $112$&$e^{\epsilon 96}205$&76800\\
\hline
$7$ & $e^{\epsilon 1}5$ & $112$&$e^{\epsilon 16}5$&76800\\
\hline
$10$ & $e^{\epsilon 2}5$ & $96$&&\\
&$e^{\epsilon 5}11$&&&\\
&$e^{\epsilon 7}15$&&$e^{\epsilon 244}511$&$124416$\\
\hline
\end{tabular}
\caption{A set of consonances in $\mathbb{Z}_{16}$, their respective counterpoint symmetries and
number of admitted successors, and their extended counterpoint symmetries when embedded
in $\mathbb{Z}_{512}$, with the corresponding number of extended admitted successors.}
\label{T:Sim}
\end{table}

\section{A possible continuous counterpoint}

The previous calculations suggest the following constructions that enable a continuous and compositionally useful counterpoint.
First, we consider the space $S^{1}\subseteq \mathbb{C}$ (which represents the continuum of intervals
modulo octave), with the action of 
the group $G=\mathbb{R}/\mathbb{Z}\ltimes \mathbb{Z}_{2}$ given by
\[
e^{t}v(x) = \begin{cases}
x\exp(2\pi it), &v=1,\\
\overline{x}\exp(2\pi it), &v=-1.\\
\end{cases}
\]

We define the set of consonances $(K/D)$ as the image of $[0,\frac{1}{2})$ under the
map $\phi:[0,1]\mapsto S^{1}:t\mapsto e^{2i\pi t}$, which musically means that we consider as consonant any interval greater or equal than the
unison but smaller than the tritone (within an octave). Apart from the identity, no element of $G$
leaves $(K/D)$ invariant, thus it is strong and its polarity is $e^{\tfrac{1}{2}}$.

Now, for counterpoint, we consider the torus $T=S^{1}\times S^{1}$, with the
first component for the cantus firmus and the second for the discantus interval. Let
$G$ act on $T$ in the following manner:
\[
e^{t}v(x,y) = (vx,e^{t}vy);
\]
this action is suggested by the fact that all the linear parts of the affine symmetries
of counterpoint intervals have no dual component.

Thus the set of consonant intervals is $(K[\epsilon]/D[\epsilon])=(S^{1}\times K/S^{1}\times D)$, the
self-complementary function for any $\xi \in T$ which fixes its tangent space is $e^{1/2}1$, and
it commutes with any element of $G'$. Also $\xi = (0,k) \in g(D[\epsilon])$ for a $g\in G'$
if and only if
\[
g = e^{t}1, \quad t\in (k,k+1/2]\quad
\text{or}\quad g = e^{t}(-1), \quad t\in [k-1/2,k).
\]

And here comes a delicate point. If we wish to preserve the idea of cardinality maximization,
it would be reasonable to ask the set of infinite admitted successors to attain certain maximum. A
possibility is to gauge these sets in terms
of the standard measure in $T$ since, for instance, the affine morphisms
\[
g=\begin{cases}
e^{k-1/2}(-1), &k\in \phi([0,1/4]),\\
e^{k-1/2}1, &k\in \phi([1/4,1/2]),
\end{cases}
\]
maximize the measure of the intersection $(gX[\epsilon])\cap X[\epsilon]$. The musical meaning
of this alternative is that the admitted successors of
consonant intervals below the minor third are all the consonant intervals above it, and
vice versa. The minor third is special, because it has any consonant interval as an admitted successor.

But, in terms of the new perspective of homology introduced by Mazzola in \cite{MR2915169}, we observe first that $T$ is homeomorphic
to $T$ itself with respect to the Kuratowski closure operator induced by the quasipolarity $e^{1/2}1$. This is so because, for in each section $x\times S^{1}$,
the self-complementary function is the antipodal morphism, thus each $x\times S^{1}$ is homeomorphic
to the projective line, which in turn is homeomorphic to $x\times S^{1}$ itself \cite[p. 58]{mC06}.
Furthermore, any $g\in G'$ which leaves $\xi$ out of
$g(X[\epsilon])$ is such that $(g(X[\epsilon]))\cap X[\epsilon]$ is homotopically equivalent to $S^{1}$,
except when such intersection is empty\footnote{By the way, this happens only with the
self-complementary function itself.}. Therefore, $H_{1}((g(X[\epsilon]))\cap X[\epsilon])=\mathbb{Z}$
is always the group of maximum rank when it satisfies the rest of the conditions
of counterpoint symmetries. This implies that, except for itself, any counterpoint interval
can be an admitted counterpoint successor, which is clearly an undesirable outcome.

\section{Some final remarks}

In the version of infinite counterpoint that maximizes measure, we arrive to
some peculiar features:
\begin{enumerate}
\item Certainly there are no culs-de-sac.
\item The only consonance that has all the other consonances as admitted successors is
the minor third.
\item All the intervals smaller than the minor third admit only larger intervals as successors,
while all those greater admit only smaller ones.
\item Although it is continuous regarding its induced quasipolarity and the cantus firmus can
be chosen to be a continuous function of time, the discantus cannot be continuous in the standard
topology.
\end{enumerate}

All of these seem to be very close to the general principles of counterpoint. Unfortunately,
this specific instance is not a natural extension of the discrete version; their relation is mainly axiomatic.
On the other hand, the restriction of the linear part
of the morphisms to $\mathbb{Z}_{2}$, although not entirely artificial, feels too limited with respect to
the original finite model.

In fact, the selected dichotomy for the continuous example is a particularly nice one that permits
a simple analysis, but by no means it is the only possible one. Considering that the general
linear parts for counterpoint symmetries can be recovered as ``windings'' of $S^{1}$, carefully
constructed infinite dichotomies could yield more complicated homology groups that make the 
algebro-topological approach far more interesting.

\bibliographystyle{plain}
\bibliography{tesis}

\end{document}